\documentclass[11pt]{amsart}
\usepackage{graphicx}
\usepackage{amssymb}

\newtheorem{thm}{Theorem}

\newtheorem{exmp}[thm]{Example}

\def\sp {\hskip .05in}
\def\bsp {\hskip .5in}
\def\bspp {\hskip .5in}

\title {A bijective proof of an unusual symmetric
group generating function}

\author{Mike Zabrocki}

\address{Mathematics and Statistics, TEL Building, York University, 4700 Keele Street, Toronto, Ontario, M3J 1P3, Canada}
\begin{document}

\begin{abstract}
For
$\sigma \in S_n$, let $D(\sigma) = \{ i : \sigma_{i} > \sigma_{i+1}\}$
denote the descent set of $\sigma$.  The length of the permutation is
the number of inversions, denoted by $inv(\sigma) = \big | \{ (i,j) : i<j, 
\sigma_i > \sigma_j\} \big |$.  Define an unusual quadratic statisitic by
$baj(\sigma) = \sum_{i \in D(\sigma)} i ( n-i)$.  We present here
a bijective proof of the identity $\sum_{{\sigma \in S_n}
\atop {\sigma(n) = k}} q^{baj(\sigma) -
inv(\sigma)} = \prod_{i=1}^{n-1} {{1-q^{i (n-i)}} \over {1-q^i}}$ 
where $k$ is a fixed integer.
\end{abstract}

\maketitle

The following identity was presented as a
special case of a Weyl group generating function
in a seminar talk  at UCSD in November 1996
by John Stembridge.  We present here a bijective proof
of the identity.

Define a statistic on the permutations on $n$ letters
\begin{eqnarray}
baj(\sigma) = \sum_{i=1}^{n-1} i (n-i) \chi ( \sigma_{i+1}
< \sigma_{i} ) \nonumber
\end{eqnarray}
 where $\chi$ is an indicator function

\begin{eqnarray} \nonumber
\chi(A) = \Bigg\{ {{1 \sp if \sp A \sp true} \atop
{0 \sp if \sp A \sp false}} 
\end{eqnarray}

The number of inversions of a permutation may be expressed as 
\begin{eqnarray} \nonumber
inv(\sigma) = \sum_{i=1}^{n-1} \sum_{j>i} \chi( \sigma_j < \sigma_i )
\end{eqnarray}

A special case of the main result presented in \cite{SW}
when the root system is $A_{n-1}$ is the following

\begin{thm}
\begin{equation} \nonumber
\sum_{\sigma \in S_n} q^{baj(\sigma) - inv(\sigma)} = 
n \prod_{i=1}^{n-1} {{1 - q^{i ( n-i )}} \over {1 - q^i}}
\end{equation}
\end{thm}

A slightly stronger statement can be made for the specialization
of this formula to this root system. For a fixed
$k \in \{ 1 \ldots n \}$ the following equation also holds

\begin{thm}
\begin{equation} \nonumber
\sum_{{\sigma \in S_n} \atop {\sigma_n = k}} q^{baj(\sigma) - inv(\sigma)} = 
\prod_{i=1}^{n-1} {{1 - q^{i ( n-i )}} \over {1 - q^i}}
\end{equation}

where $k$ is an integer between $1$ and $n$.
\end{thm}


Notice that the right hand side of this equation can be expressed
as a product of sums by the formula

\begin{eqnarray} \nonumber
\prod_{i=1}^{n-1} {{1 - q^{i ( n-i )}} \over {1 - q^i}}
&=& \prod_{i=1}^{n-1} {{1 - q^{i ( n-i )}} \over {1 - q^{n-i}}}
= \prod_{i=1}^{n-1} \sum_{r_i = 0}^{i-1} q^{(n-i) r_i} \\
&=& \sum_{{{(r_1, r_2, \ldots, r_{n-1})} \atop {0 \leq r_i < i}}}
q^{\sum_{i=1}^{n-1} (n-i) r_i} \nonumber
\end{eqnarray}

The object of this proof will be to find a
bijection from the permutations, $\sigma$, of $\{1 \ldots n \}$ 
with $\sigma_n = k$ ($k$ fixed) to 
sequences of integers $(r_1, r_2, \ldots r_{n-1})$ with the
additional property that
\begin{eqnarray} \nonumber
baj(\sigma) - inv(\sigma) = \sum_{i=1}^{n-1} (n-i) r_i
\end{eqnarray}

Let $\sigma$ be a permutation of $\{1 \ldots n\}$.
$\sigma$ can be represented by a sequence of
integers $(v_1, v_2, \ldots, v_n)$ where $v_i$ is 
$1 \leq v_i \leq i$ and is given by $v_i =
\sum_{j \leq i} \chi( \sigma_j \leq \sigma_i ) $.

Given such a sequence of $v_i$, it is possible to recover the
permutation that it corresponds to
by first constructing $\sigma' \in S_{n-1}$ for the sequence
$(v_1, v_2, \ldots, v_{n-1})$ and then defining the permutation
$\sigma \in S_n$
by $\sigma_n = v_n$, 
$\sigma_i = \sigma_i'$ if $\sigma_i' < v_n$,
and $\sigma_i = \sigma_i' +1$
if $\sigma_i \geq v_n$.  This construction gives that the number
of $i \in \{ 1 \ldots n \}$ 
such that $\sigma_i \leq \sigma_n$ is $v_n$.  This quantity does not
change by building a larger permutation in the same manner.

\begin{exmp}

Say $v = (1,1,3,1,2,5,1)$ \\
$\bspp v_1 = 1 \bsp \sigma^{(1)} = 1$ \\
$\bspp v_2 = 1 \bsp \sigma^{(2)} = 21$ \\
$\bspp v_3 = 3 \bsp \sigma^{(3)} = 213$ \\
$\bspp v_4 = 1 \bsp \sigma^{(4)} = 3241$ \\
$\bspp v_5 = 2 \bsp \sigma^{(5)} = 43512$ \\
$\bspp v_6 = 5 \bsp \sigma^{(6)} = 436125$ \\
$\bspp v_7 = 1 \bsp \sigma^{(7)} = 5472361 = \sigma$
\end{exmp}

Notice that $\sigma_{i+1} < \sigma_i$ if and only if $v_{i+1} \leq
v_i$.  This is because if $v_{i+1} \leq v_i$ then $\sigma^{(i+1)}$ will
have a descent in the $i^{th}$ position, and this descent will remain
for all $\sigma^{(k)}$ with $k>i$
(and in particular $\sigma^{(n)} = \sigma$).

Define the bijection from permutations $\sigma$ with $\sigma_n = k$ by
first computing the sequence of $v_i$ and then setting
$r_i = i \chi(v_{i+1} \leq v_i) + v_{i+1} - v_i -1$.

This defines a map from such permutations to sequences of integers
$(r_1, r_2, \ldots, r_{n-1})$ with $0 \leq r_i < i$.
Note that if $v_i \geq v_{i+1}$ then $i-1 \geq v_i - v_{i+1} \geq 0$
so that $0 \leq r_i = i-1 -(v_i - v_{i+1}) \leq i-1$. If
$v_i < v_{i+1}$ then $0 < v_{i+1} - v_i \leq i$, hence
$0 \leq r_i = v_{i+1} - v_i - 1 \leq i-1$.

Given a sequence $(r_1, r_2, \ldots, r_{n-1})$ that is
the image of some permutation and assume
that the values of
$v_{i+1}, \ldots, v_{n}$ are known, then $v_i - i \chi(v_{i+1} \leq v_i) =
v_{i+1} - r_i + 1$.  If the right hand side of the equation less than
or equal to $0$ then it must be that $\chi(v_{i+1} \leq v_i) =1$
and so $v_i = i + v_{i+1} - r_i  + 1$.  Otherwise
$\chi(v_{i+1} \leq v_i) =0$ and then $v_i = v_{i+1} - r_i + 1$.
The whole sequence of $v_i$ can be recovered, and thus, the original
permutation also.

This map is 1-1 since it is possible to recover
the permutation if 
the sequence $(r_1, r_2, \ldots, r_{n-1})$ is given and
the value of $v_n = k$ is known.  There are
the same number of permutations with the last element fixed as there
are such sequences of numbers, hence this map is a bijection.

It remains to show the result
\begin{eqnarray} \nonumber
baj(\sigma) - inv(\sigma) = \sum_{i=1}^{n-1} (n-i) r_i
\end{eqnarray}

Note that $v_i$ can be expressed by the formula
$v_i = \sum_{j \leq i} \chi( \sigma_j \leq \sigma_i )$. 
Because
$\big({{n+1} \atop {2}} \big) =
\sum_{i=1}^{n} \sum_{j \leq i} (\chi( \sigma_j \leq \sigma_i ) +
\chi(\sigma_j > \sigma_i) ) = \sum_{i=1}^{n} v_i + inv(\sigma)$,
the number of inversions of the permutation $\sigma$ is given
by the formula
\begin{eqnarray} \nonumber
inv(\sigma) = \bigg({{n+1} \atop {2}} \bigg) - \sum_{i=1}^n v_i
\end{eqnarray}

The statistic $baj$ can be given in terms of the $v_i$'s because of
the remark that $v_{i+1} \leq v_i$ if and only if $\sigma_{i+1} < \sigma_i$.
\begin{eqnarray} \nonumber
baj( \sigma) = \sum_{i=1}^{n-1} i ( n-i) \chi(v_{i+1} \leq v_i)
\end{eqnarray}

Thus,
\begin{eqnarray} \nonumber
\sum_{i=1}^{n-1} r_i(n-i) &=& \sum_{i=1}^{n-1} (i \chi(
v_{i+1} 
\leq v_i) + v_{i+1} - v_i -1)(n-i)  \nonumber \\
&=& \sum_{i=1}^{n-1} i(n-i) \chi( v_{i+1} \leq v_i) + \sum_{i=1}^{n-1}
v_{i+1} (n-i) - \sum_{i=1}^{n-1} (v_i+1)(n-i) 
\nonumber \\
&=& \sum_{i=1}^{n-1} i(n-i) \chi( v_{i+1} \leq v_i) +
\sum_{i=2}^{n} v_{i} (n-i+1) - \sum_{i=1}^{n-1} v_i (n-i) -
\sum_{i=1}^{n-1} i  \nonumber \\
&=& \sum_{i=1}^{n-1} i(n-i) \chi( v_{i+1}
\leq v_i) + v_n +
\sum_{i=2}^{n-1} v_{i} -  v_1 (n-1) -  \bigg({{n} \atop {2}} \bigg) 
 \nonumber \\
&=& \sum_{i=1}^{n-1} i(n-i) \chi( v_{i+1} \leq v_i) +
\sum_{i=1}^{n} v_{i} -  v_1 n -  \bigg({{n} \atop {2}} \bigg) 
 \nonumber \\
&=& \sum_{i=1}^{n-1} i ( n-i) \chi(v_{i+1} \leq v_i) +
\sum_{i=1}^n v_i - \bigg({{n+1} \atop {2}} \bigg) 
 \nonumber \\
&=& \, baj(\sigma) - inv(\sigma)  \nonumber
\end{eqnarray}

This shows the last property of the bijection and hence the theorem.

\end{document}